\documentclass[11pt,reqno]{amsart}

\usepackage{amsmath,amsthm,amscd}
\usepackage{amssymb}
\usepackage{txfonts}
\usepackage{euscript}
\usepackage[all]{xy}
\usepackage{epic,eepic}

\numberwithin{equation}{subsection}

\newcommand{\sgsp}{\renewcommand{\baselinestretch}{1}\tiny\normalsize}
\newcommand{\sqsp}{\renewcommand{\baselinestretch}{1.2}\tiny\normalsize}

\newtheorem{thm}[subsection]{Theorem}
\newtheorem{lemma}[subsection]{Lemma}
\newtheorem{prop}[subsection]{Proposition}

\newcommand{\UHLie}{U_{HLie}}
\newcommand{\UHLeib}{U_{HLeib}}
\newcommand{\FHNA}{F_{HNAs}}
\newcommand{\FHAs}{F_{HAs}}

\newcommand{\HomMod}{\mathbf{HomMod}}
\newcommand{\HomLie}{\mathbf{HomLie}}
\newcommand{\HomLeib}{\mathbf{HomLeib}}
\newcommand{\HomAs}{\mathbf{HomAs}}
\newcommand{\HomDi}{\mathbf{HomDi}}
\newcommand{\HomNonAs}{\mathbf{HomNonAs}}

\newcommand{\Fbar}{\overline{F}}
\newcommand{\Twt}{T^{wt}}
\newcommand{\Tdi}{T^{di}}

\newcommand{\lprod}{\dashv}  % -|
\newcommand{\rprod}{\vdash}  % |-

\DeclareMathOperator{\Id}{Id}
\DeclareMathOperator{\Hom}{Hom}
\DeclareMathOperator{\End}{End}

\begin{document}
\title{Enveloping algebras of Hom-Lie algebras}
\author{Donald Yau}

\begin{abstract}
Enveloping algebras of Hom-Lie and Hom-Leibniz algebras are constructed.
\end{abstract}


\email{dyau@math.ohio-state.edu}
\address{Department of Mathematics, The Ohio State University at Newark, 1179 University Drive, Newark, OH 43055, USA}

%\date{\today}
\maketitle

\sqsp

%%================%%
%%================%%
%%                %%
%%  Introduction  %%
%%                %%
%%================%%
%%================%%

\section{Introduction}
\label{sec:intro}

% Hom-Lie algebra

A Hom-Lie algebra is a triple $(L, \lbrack -,- \rbrack, \alpha)$, where $\alpha$ is a linear self-map, in which the skew-symmetric bracket satisfies an $\alpha$-twisted variant of the Jacobi identity, called the Hom-Jacobi identity.  When $\alpha$ is the identity map, the Hom-Jacobi identity reduces to the usual Jacobi identity, and $L$ is a Lie algebra.  Hom-Lie algebras and related algebras were introduced in \cite{hls} to construct deformations of the Witt algebra, which is the Lie algebra of derivations on the Laurent polynomial algebra $\mathbf{C}\lbrack z^{\pm 1} \rbrack$.

% Hom-associative algebra

An elementary but important property of Lie algebras is that each associative algebra $A$ gives rise to a Lie algebra $Lie(A)$ via the commutator bracket.  In \cite{ms}, Makhlouf and Silvestrov introduced the notion of a Hom-associative algebra $(A, \mu, \alpha)$, in which the binary operation $\mu$ satisfies an $\alpha$-twisted version of associativity.  Hom-associative algebras play the role of associative algebras in the Hom-Lie setting.  In other words, a Hom-associative algebra $A$ gives rise to a Hom-Lie algebra $HLie(A)$ via the commutator bracket.

% Enveloping algebras for Hom-Lie algebra

The first main purpose of this paper is to construct the enveloping Hom-associative algebra $\UHLie(L)$ of a Hom-Lie algebra $L$.  In other words, $\UHLie$ is the left adjoint functor of $HLie$.  This is analogous to the fact that the functor $Lie$ admits a left adjoint $U$, the enveloping algebra functor.  The construction of $\UHLie(L)$ makes use of the combinatorial objects of weighted binary trees, i.e., planar binary trees in which the internal vertices are equipped with weights of non-negative integers.

% Hom-Leibniz and Hom-dialgebras

The second main purpose of this paper is to construct the counterparts of the functors $HLie$ and $\UHLie$ for Hom-Leibniz algebras.  Leibniz algebras  (also known as \emph{right Loday algebras}) \cite{loday0,loday1,loday,lp,lp2} are non-skew-symmetric versions of Lie algebras in which the bracket satisfies a variant of the Jacobi identity.  In particular, Lie algebras are examples of Leibniz algebras.  In the Leibniz setting, the objects that play the role of associative algebras are called dialgebras, which were introduced by Loday in \cite{loday}.  Dialgebras have two associative binary operations that satisfy three additional associative-type axioms.  Leibniz algebras were extended to Hom-Leibniz algebras in \cite{ms}.  Extending some of the work of Loday \cite{loday}, we will introduce \emph{Hom-dialgebras} and construct the adjoint pair $(\UHLeib, HLeib)$ of functors.  In a Hom-dialgebra, there are two binary operations that satisfy five $\alpha$-twisted associative-type axioms.

Each Hom-Lie algebra can be thought of as a Hom-Leibniz algebra.  Likewise, to every Hom-associative algebra is associated a Hom-dialgebra in which both binary operations are equal to the original one.  The functors described above give rise to the following diagram of categories and (adjoint) functors.
   \[
   \UseTips
   \newdir{ >}{!/-5pt/\dir{>}}
   \xymatrix{
   \HomLie \ar@<.5ex>[rr]^{\UHLie} \ar@{^{(}->}[d]_{\iota} & &  \HomAs \ar@<.5ex>[ll]^{HLie} \ar@{^{(}->}[d]_{\iota} \\
   \HomLeib \ar@<.5ex>[rr]^{\UHLeib} & & \HomDi \ar@<.5ex>[ll]^{HLeib}
   }
   \]
Moreover, the commutativity, $HLeib \circ \iota = \iota \circ HLie$, holds.

% Motivation: purpose of UHLie and UHleib

This paper is the first part of a bigger project to study (co)homology theories of the various Hom-algebras.  In papers under preparation, the author aims to construct:
\begin{itemize}
\item Hochschild-type (co)homology and its corresponding cyclic (co)homology for Hom-associative algebras,
\item an analogue of the Chevalley-Eilenberg Lie algebra (co)homology for Hom-Lie algebras, and
\item an analogue of Loday's Leibniz algebra (co)homology for Hom-Leibniz algebras.
\end{itemize}
From the point-of-view of homological algebra, it is often desirable to interpret homology and cohomology in terms of resolutions and the derived functors $Tor$ and $Ext$, respectively.  In the classical case of Lie algebra (co)homology, this requires the enveloping algebra functor $U$, since $H_*^{Lie}(L,-) \cong Tor^{UL}_*(\mathbb{K},-)$ and $H^*_{Lie}(L,-) \cong Ext^*_{UL}(\mathbb{K},-)$ for a Lie algebra $L$, where $\mathbb{K}$ is the ground field.  It is reasonable to expect that our enveloping algebra functors $\UHLie$ and $\UHLeib$, or slight variations of these functors, will play similar roles for Hom-Lie and Hom-Leibniz (co)homology, respectively.

A major reason to study Hom-Lie algebra cohomology is to provide a proper context for the \emph{Hom-Lie algebra extensions} constructed in \cite[Section 2.4]{hls}.  After constructing a certain $q$-deformation of the Witt algebra, Hartwig, Larsson, and Silvestrov used the machinery of Hom-Lie algebra extensions to construct a corresponding deformation of the Virasoro algebra \cite[Section 4]{hls}.  In the classical case of Lie algebras, equivalence classes of extensions are classified by the cohomology module $H^2_{Lie}$.  The author expects that Hom-Lie algebra extensions will admit a similar interpretation in terms of the second Hom-Lie algebra cohomology module.

A important result in the theory of Lie algebra homology is the Loday-Quillen Theorem \cite{lq} relating Lie algebra and cyclic homology.  In another paper, the author hopes to extend the Loday-Quillen Theorem to the Hom-algebra setting.  Recall that the Loday-Quillen Theorem states that, for an associative algebra $A$, $H^{Lie}_*(gl(A),\mathbb{K})$ is isomorphic to the graded symmetric algebra of $HC_{*-1}(A)$.  It is not hard to check that, if $A$ is a Hom-associative algebra, then so is $gl(A)$.  Therefore, using the commutator bracket of Makhlouf and Silvestrov \cite{ms}, $gl(A)$ can also be regarded as a Hom-Lie algebra.  There is an analogue of the Loday-Quillen Theorem for Leibniz algebra and Hochschild homology due to Loday \cite[Theorem 10.6.5]{loday0}, which the author also hopes to extend to the Hom-algebra case.

% Organization

\subsection{Organization}

The rest of this paper is organized as follows.

The next section contains preliminary materials on binary trees.  In Section \ref{sec:HomMod}, free Hom-nonassociative algebras of Hom-modules are constructed (Theorem \ref{thm:FreeHomNonAs}).  This leads to the construction of the enveloping Hom-associative algebra functor $\UHLie$ in Section \ref{sec:env} (Theorem \ref{thm:env}).

In Section \ref{sec:HomDi}, Hom-dialgebras are introduced together with several classes of examples.  It is then observed that Hom-dialgebras give rise to Hom-Leibniz algebras via a version of the commutator bracket (Proposition \ref{prop:HomDi2HomLeib}).  The enveloping Hom-dialgebra functor $\UHLeib$ for Hom-Leibniz algebras are constructed in Section \ref{sec:env Leibniz} (Theorem \ref{thm:envLeib}).

%In Section \ref{sec:free}, free $G$-Hom-associative, free Hom-Lie, and free Hom-Leibniz algebras as well as free Hom-dialgebras are constructed.

%%%%%%%%%%%%%%%%%%%%%%%%%%%%%%%%%%%%%%%
%%%%%%%%%%%%%%%%%%%%%%%%%%%%%%%%%%%%%%%
%%
%%     Binary Trees
%%
%%%%%%%%%%%%%%%%%%%%%%%%%%%%%%%%%%%%%%%
%%%%%%%%%%%%%%%%%%%%%%%%%%%%%%%%%%%%%%%
\section{Preliminaries on binary trees}
\label{sec:trees}

The purpose of this section is to collect some basic facts about binary trees that are needed for the construction of the enveloping algebra functors in later sections.  Sections \ref{subsec:trees} and \ref{subsec:grafting} below follow the discussion in \cite[Appendix A1]{loday} but with slightly different notation.

%%%%%%%%%%%%%%%%%%%%%%%%%%%%%%%%
\subsection{Planar binary trees}
\label{subsec:trees}

For $n \geq 1$, let $T_n$ denote the set of planar binary trees with $n$ leaves and one root.  Below are the first four sets $T_n$.
   \[
   T_1 = \left\lbrace \begin{picture}(2,10)      % T_1
                 \drawline(1,0)(1,10)
                 \end{picture}\,
         \right\rbrace,
   T_2 = \left\lbrace \begin{picture}(16,14)     % T_2
                 \drawline(8,0)(8,6)(0,14)
                 \drawline(8,6)(16,14)
                 \put(8,6){\circle*{2}}
                 \end{picture}
         \right\rbrace,
   T_3 = \left\lbrace \begin{picture}(16,14)     % T_3[1]
                 \drawline(8,0)(8,6)(0,14)
                 \drawline(8,6)(16,14)
                 \drawline(12,10)(8,14)
                 \put(8,6){\circle*{2}}
                 \put(12,10){\circle*{2}}
                 \end{picture},\,
                 \begin{picture}(16,14)     % T_3[2]
                 \drawline(8,0)(8,6)(0,14)
                 \drawline(4,10)(8,14)
                 \drawline(8,6)(16,14)
                 \put(8,6){\circle*{2}}
                 \put(4,10){\circle*{2}}
                 \end{picture}
         \right\rbrace,
   T_4 = \left\lbrace \begin{picture}(18,15)(0,-4)    % T_4[1]
                 \drawline(9,-4)(9,2)(0,11)
                 \drawline(9,2)(18,11)
                 \drawline(12,5)(6,11)
                 \drawline(15,8)(12,11)
                 \put(9,2){\circle*{2}}
                 \put(12,5){\circle*{2}}
                 \put(15,8){\circle*{2}}
                 \end{picture},\,
                 \begin{picture}(18,15)(0,-4)    % T_4[2]
                 \drawline(9,-4)(9,2)(0,11)
                 \drawline(9,2)(18,11)
                 \drawline(12,5)(6,11)
                 \drawline(9,8)(12,11)
                 \put(9,2){\circle*{2}}
                 \put(12,5){\circle*{2}}
                 \put(9,8){\circle*{2}}
                 \end{picture},\,
                 \begin{picture}(18,15)(0,-4)    % T_4[3]
                 \drawline(9,-4)(9,2)(0,11)
                 \drawline(3,8)(6,11)
                 \drawline(9,2)(18,11)
                 \drawline(15,8)(12,11)
                 \put(9,2){\circle*{2}}
                 \put(3,8){\circle*{2}}
                 \put(15,8){\circle*{2}}
                 \end{picture},\,
                 \begin{picture}(18,15)(0,-4)    % T_4[4]
                 \drawline(9,-4)(9,2)(0,11)
                 \drawline(6,5)(12,11)
                 \drawline(9,8)(6,11)
                 \drawline(9,2)(18,11)
                 \put(9,2){\circle*{2}}
                 \put(6,5){\circle*{2}}
                 \put(9,8){\circle*{2}}
                 \end{picture},\,
                 \begin{picture}(18,15)(0,-4)    % T_4[5]
                 \drawline(9,-4)(9,2)(0,11)
                 \drawline(3,8)(6,11)
                 \drawline(6,5)(12,11)
                 \drawline(9,2)(18,11)
                 \put(9,2){\circle*{2}}
                 \put(6,5){\circle*{2}}
                 \put(3,8){\circle*{2}}
                 \end{picture}
          \right\rbrace.
   \]
Each dot represents an internal vertex.  From now on, an element of $T_n$ will simply be called an \emph{$n$-tree}.  In each $n$-tree, there are $(n-1)$ internal vertices, the lowest one of which is connected to the root.  Note that the cardinality of the set $T_{n+1}$ is the $n$th Catalan number $C_n = (2n)!/n!(n+1)!$.

%%%%%%%%%%%%%%%%%%%%%%%%%%%%%%
\subsection{Grafting of trees}
\label{subsec:grafting}

Let $\psi \in T_n$ and $\varphi \in T_m$ be two trees.  Define an $(n+m)$-tree $\psi \vee \varphi$, called the \emph{grafting} of $\psi$ and $\varphi$, by joining the roots of $\psi$ and $\varphi$ together, which forms the new lowest internal vertex that is connected to the new root.  Pictorially, we have
   \[
   \psi \vee \varphi \,=\,
                 \begin{picture}(20,18)(-2,0)
                 \drawline(8,0)(8,6)(0,14)
                 \drawline(8,6)(16,14)
                 \put(8,6){\circle*{2}}
                 \put(-1,18){\makebox(0,0){$\psi$}}
                 \put(19,18){\makebox(0,0){$\varphi$}}
                 \end{picture}.
   \]
Note that grafting is a nonassociative operation.

Conversely, by cutting the two upward branches from the lowest internal vertex, each $n$-tree $\psi$ can be uniquely represented as the grafting of two trees, say, $\psi_1 \in T_p$ and $\psi_2 \in T_q$, where $p + q = n$.  By iterating the grafting operation, one can show by a simple induction argument that every $n$-tree $(n \geq 2)$ can be obtained as an iterated grafting of $n$ copies of the $1$-tree.

%%%%%%%%%%%%%%%%%%%%%%%%%%%%%
\subsection{Weighted trees}
\label{subsec:weighted trees}

By a \emph{weighted $n$-tree}, we mean a pair $\tau = (\psi, w)$, in which:
   \begin{enumerate}
   \item $\psi \in T_n$ is an $n$-tree and
   \item $w$ is a function from the set of internal vertices of $\psi$ to the set $\mathbb{Z}_{\geq 0}$ of non-negative integers.
   \end{enumerate}
If $v$ is an internal vertex of $\psi$, then we call $w(v)$ the \emph{weight} of $v$.  The $n$-tree $\psi$ is called the \emph{underlying $n$-tree} of $\tau$, and $w$ is called the \emph{weight function} of $\tau$.  Let $\Twt_n$ denote the set of all weighted $n$-trees.  Since the $1$-tree has no internal vertex, we have that $T_1 = \Twt_1$.

We can picture a weighted $n$-tree $\tau = (\psi, w)$ by drawing the underlying $n$-tree $\psi$ and putting the weight $w(v)$ next to each internal vertex $v$.  For example, here is a weighted $4$-tree,
   \begin{equation}
   \label{eq:4tree}
   \tau \,=\,
   \setlength{\unitlength}{1mm}
   \begin{picture}(18,15)(0,-1)    % T_4[2]
                 \drawline(9,-1)(9,2)(0,11)
                 \drawline(9,2)(18,11)
                 \drawline(12,5)(6,11)
                 \drawline(9,8)(12,11)
                 \put(9,2){\circle*{.8}}
                 \put(12,5){\circle*{.8}}
                 \put(9,8){\circle*{.8}}
                 \put(7.1,7.4){\makebox(0,0){$5$}}
                 \put(14.5,4.5){\makebox(0,0){$2$}}
                 \put(7.2,1.3){\makebox(0,0){$7$}}
   \end{picture}.
   \end{equation}

%%%%%%%%%%%%%%%%%%%%%%%%%%%%%%%%%%%%%%%%%
\subsection{Grafting of weighted trees}
\label{subsec:graft wt}

Let $\tau = (\psi, w) \in \Twt_p$ and $\tau^\prime = (\psi^\prime, w^\prime) \in \Twt_q$ be two weighted trees.  Define their \emph{grafting} to be the weighted $(p+q)$-tree
   \[
   \tau \vee \tau^\prime \,\buildrel \text{def} \over=\,
   (\psi \vee \psi^\prime, \omega)
   \]
with underlying tree $\psi \vee \psi^\prime$.  The weight function is given by
   \[
   \omega(v) \,=\,
   \begin{cases}
   w(v) & \text{if $v$ is an internal vertex of $\psi$}, \\
   w^\prime(v) & \text{if $v$ is an internal vertex of $\psi^\prime$}, \\
   0 & \text{if $v$ is the lowest internal vertex of $\psi \vee \psi^\prime$}.
   \end{cases}
   \]
The grafting can be pictured as
   \[
   \tau \vee \tau^\prime \,=\,
   \setlength{\unitlength}{.5mm}
                 \begin{picture}(20,18)(-2,0)
                 \drawline(8,0)(8,6)(0,14)
                 \drawline(8,6)(16,14)
                 \put(8,6){\circle*{2}}
                 \put(13,5){\makebox(0,0){$0$}}
                 \put(-2,15.6){$\tau$}
                 \put(18.5,18){\makebox(0,0){$\tau^\prime$}}
                 \end{picture}.
   \]

%%%%%%%%%%%%%%%%%%%%%%%%%%%%%%%
\subsection{The $+m$ operation}
\label{subsec:1}

Let $\tau = (\psi, w) \in \Twt_p$ be a weighted $p$-tree, and let $m$ be a non-negative integer.  Define a new weighted $p$-tree
   \[
   \tau\lbrack m \rbrack \, \buildrel\text{def}\over =\,
   (\psi, w\lbrack m \rbrack)
   \]
with the same underlying $p$-tree $\psi$.  The weight function is given by
   \[
   w\lbrack m \rbrack(v) \,=\,
   \begin{cases}
   w(v) & \text{ if $v$ is not the lowest internal vertex of $\psi$}, \\
   w(v) + m & \text{ if $v$ is the lowest internal vertex of $\psi$}.
   \end{cases}
   \]
Pictorially, if
   \[
   \tau \,=\,
   \setlength{\unitlength}{.5mm}
                 \begin{picture}(20,18)(-2,0)
                 \drawline(8,0)(8,6)(0,14)
                 \drawline(8,6)(16,14)
                 \put(8,6){\circle*{2}}
                 \put(13,5){\makebox(0,0){$r$}}
                 \put(-1,18){\makebox(0,0){$\tau_1$}}
                 \put(19,18){\makebox(0,0){$\tau_2$}}
                 \end{picture},
   \]
then
   \[
   \tau\lbrack m \rbrack \,=\,
   \setlength{\unitlength}{.5mm}
                 \begin{picture}(30,18)(-2,0)
                 \drawline(8,0)(8,6)(0,14)
                 \drawline(8,6)(16,14)
                 \put(8,6){\circle*{2}}
                 \put(10,3){$r+m$}
                 \put(-1,18){\makebox(0,0){$\tau_1$}}
                 \put(19,18){\makebox(0,0){$\tau_2$}}
                 \end{picture}.
   \]

By cutting the two upward branches from the lowest internal vertex, every weighted $n$-tree $\tau$ can be written uniquely as
   \begin{equation}
   \label{eq:tau}
   \tau \,=\, (\tau_1 \vee \tau_2) \lbrack r \rbrack,
   \end{equation}
where $\tau_1 \in \Twt_p$, $\tau_2 \in \Twt_q$ with $p + q = n$, and $r$ is the weight of the lowest internal vertex of $\tau$.  The same process can be applied to $\tau_1$ and $\tau_2$, and so on.  In particular, every weighted $n$-tree for $n \geq 2$ can be obtained from $n$ copies of the $1$-tree by iterating the operations $\vee$ and $\lbrack r \rbrack$ $(r \geq 0)$.  For example, denoting the $1$-tree by $i$, the weighted $4$-tree in \eqref{eq:4tree} can be written as
   \[
   \setlength{\unitlength}{1mm}
   \begin{picture}(18,15)(0,-1)    % T_4[2]
                 \drawline(9,-1)(9,2)(0,11)
                 \drawline(9,2)(18,11)
                 \drawline(12,5)(6,11)
                 \drawline(9,8)(12,11)
                 \put(9,2){\circle*{.8}}
                 \put(12,5){\circle*{.8}}
                 \put(9,8){\circle*{.8}}
                 \put(7.1,7.4){\makebox(0,0){$5$}}
                 \put(14.5,4.5){\makebox(0,0){$2$}}
                 \put(7.2,1.3){\makebox(0,0){$7$}}
   \end{picture}
   \,=\,
   \left \lbrace
   i \vee \left((i \vee i) \lbrack 5 \rbrack \vee i\right) \lbrack 2 \rbrack
   \right \rbrace
   \lbrack 7 \rbrack.
   \]

%%%%%%%%%%%%%%%%%%%%%%%%%%%%%%%%%%%%%%%%%%%%%%%%%%%%%
%%%%%%%%%%%%%%%%%%%%%%%%%%%%%%%%%%%%%%%%%%%%%%%%%%%%%
%%
%%     Hom-modules
%%
%%%%%%%%%%%%%%%%%%%%%%%%%%%%%%%%%%%%%%%%%%%%%%%%%%%%%
%%%%%%%%%%%%%%%%%%%%%%%%%%%%%%%%%%%%%%%%%%%%%%%%%%%%%
\section{Hom-modules and Hom-nonassociative algebras}
\label{sec:HomMod}

The purpose of this section is to construct the free Hom-nonassociative algebra functor on which the enveloping algebra functors are based.

Throughout the rest of this paper, let $\mathbb{K}$ denote a field of characteristic $0$. Unless otherwise specified, modules, $\otimes$, $\Hom$, and $\End$ (linear endomorphisms) are all meant over $\mathbb{K}$.

%%%%%%%%%%%%%%%%%%%%%%%%
\subsection{Hom-modules}
\label{subsec:Hom-mod}

By a \emph{Hom-module}, we mean a pair $(V, \alpha)$ consisting of:
   \begin{enumerate}
   \item a module $V$ and
   \item a linear endomorphism $\alpha \in \End(V)$.
   \end{enumerate}
A morphism $(V, \alpha) \to (V^\prime, \alpha^\prime)$ of Hom-modules is a linear map $f \colon V \to V^\prime$ such that $\alpha^\prime \circ f = f \circ \alpha$.  The category of Hom-modules is denoted by $\HomMod$.

%%%%%%%%%%%%%%%%%%%%%%%%%%%%%%%%%%%%%%%%
\subsection{Hom-nonassociative algebras}
\label{subsec:HomNA}

By a \emph{Hom-nonassociative algebra}, we mean a triple $(A, \mu, \alpha)$ in which:
   \begin{enumerate}
   \item $A$ is a module,
   \item $\mu \colon A^{\otimes 2} \to A$ is a bilinear map, and
   \item $\alpha \in \End(A)$.
   \end{enumerate}
A morphism $f \colon (A, \mu, \alpha) \to (A^\prime, \mu^\prime, \alpha^\prime)$ is a linear map $f \colon A \to A^\prime$ such that $\alpha^\prime \circ f = f \circ \alpha$ and $f \circ \mu = \mu^\prime \circ f^{\otimes 2}$.  The category of Hom-nonassociative algebras is denoted by $\HomNonAs$.

%%%%%%%%%%%%%%%%%%%%%%%%%%%%%%%%%%%%
\subsection{Parenthesized monomials}
\label{subsec:parenthesized}

In a Hom-nonassociative algebra $(A, \mu, \alpha)$, we will often abbreviate $\mu(x, y)$ to $xy$ for $x, y \in A$.  In general, given elements $x_1, \ldots , x_n \in A$, there are $\# T_n = C_{n-1}$ ways to parenthesize the monomial $x_1 \cdots x_n$ to obtain an element in $A$.  Indeed, given an $n$-tree $\psi \in T_n$, one can label the $n$ leaves of $\psi$ from left to right as $x_1, \ldots , x_n$.  Starting from the top, at each internal vertex $v$ of $\psi$, we multiply the two elements represented by the two upward branches connected to $v$.  For example, for the five $4$-trees in $T_4$ as displayed in Section \ref{subsec:trees}, the corresponding parenthesized monomials of $x_1x_2x_3x_4$ are
   \[
   x_1(x_2(x_3x_4)), ~
   x_1((x_2x_3)x_4), ~
   (x_1x_2)(x_3x_4), ~
   (x_1(x_2x_3))x_4, ~
   ((x_1x_2)x_3)x_4.
   \]
Conversely, every parenthesized monomial $x_1 \cdots x_n$ corresponds to an $n$-tree.

%%%%%%%%%%%%%%%%%%%%%
\subsection{Ideals}
\label{subsec:ideals}

Let $(A, \mu, \alpha)$ be a Hom-nonassociative algebra, and let $S \subset A$ be a non-empty subset of elements of $A$.  Then the \emph{two-sided ideal} $\langle S \rangle$ generated by $S$ is the smallest sub-$\mathbb{K}$-module of $A$ containing $S$ that is closed under $\mu$ (but is not necessarily closed under $\alpha$).  The two-sided ideal generated by $S$ always exists and can be constructed as the sub-$\mathbb{K}$-module of $A$ spanned by all the parenthesized monomials $x_1 \cdots x_n$ in $A$ with $n \geq 1$ such that at least one $x_j$ lies in $S$.  This notion of two-sided ideals will be used below (first in Section \ref{subsec:env}) in the constructions of the enveloping Hom-algebras and the free Hom-associative algebras (Section \ref{subsec:freeHomAs}).

It should be noted that a \emph{two-sided ideal} as defined in the previous paragraph is \emph{not} an ideal in the category $\HomMod$, since it is not necessarily closed under $\alpha$.

%%%%%%%%%%%%%%%%%%%%%%%%%%%%%%%%%%%%%%%%%%
\subsection{Products using weighted trees}
\label{subsec:products wt}

Each weighted $n$-tree provides a way to multiply $n$ elements in a Hom-nonassociative algebra $(A, \mu, \alpha)$.  More precisely, we define maps
   \begin{equation}
   \label{eq:product}
   \mathbb{K}\lbrack \Twt_n \rbrack \otimes A^{\otimes n} ~ \to ~ A, \quad (\tau; x_1 \otimes \cdots \otimes x_n)  ~ \mapsto ~ (x_1\cdots x_n)_\tau
   \end{equation}
inductively via the rules:
   \begin{enumerate}
   \item $(x)_i = x$ for $x \in A$, where $i$ denotes the $1$-tree.
   \item If $\tau = (\tau_1 \vee \tau_2) \lbrack r \rbrack$ as in \eqref{eq:tau}, then
      \[
      (x_1 \cdots x_n)_\tau \,=\,
      \alpha^r\left((x_1 \cdots x_p)_{\tau_1} (x_{p+1} \cdots x_{p+q})_{\tau_2}\right),
      \]
   \end{enumerate}
where $\alpha^r = \alpha \circ \cdots \circ \alpha$ ($r$ times).  This is a generalization of the parenthesized monomials discussed above.  For example, if $\tau$ is the weighted $4$-tree in \eqref{eq:4tree}, then
   \[
   (\tau; x_1 \otimes \cdots \otimes x_4) \,=\,
   \alpha^7\left(x_1\left(\alpha^2\left((\alpha^5(x_2 x_3))x_4\right)\right)\right)
   \]
for $x_1, \ldots , x_4 \in A$.  Note that
   \[
   (x_1 \cdots x_n)_{\tau\lbrack m \rbrack}
    = \alpha^m\left((x_1 \cdots x_n)_\tau \right), \,
   (x_1 \cdots x_p)_\tau (x_{p+1} \cdots x_{p+q})_\sigma
    = (x_1 \cdots x_{p+q})_{\tau \vee \sigma}.
   \]

%%%%%%%%%%%%%%%%%%%%%%%%%%%%%%%%%%%%%%%%%%%%%
\subsection{Free Hom-nonassociative algebras}
\label{subsec:free HomNA}

Let $E \colon \HomNonAs \to \HomMod$ be the forgetful functor that forgets about the binary operations.

\begin{thm}
\label{thm:FreeHomNonAs}
The functor $E$ admits a left adjoint $\FHNA \colon \HomMod \to \HomNonAs$ defined as
   \[
   \FHNA(V) \,=\, \bigoplus_{n \,\geq\, 1} \bigoplus_{\tau \,\in\, \Twt_n} \, V^{\otimes n}_\tau,
   \]
for $(V, \alpha_V) \in \HomMod$, where $V^{\otimes n}_\tau$ is a copy of $V^{\otimes n}$.
\end{thm}

\begin{proof}
First we equip $\FHNA(V)$ with the structure of a Hom-nonassociative algebra.  For elements $x_1, \ldots , x_n \in V$, a generator $(x_1 \otimes \cdots \otimes x_n)_\tau \in V^{\otimes n}_\tau$ will be abbreviated to $(x_{1,n})_\tau$.  The binary operation
   \[
   \mu_F \colon \FHNA(V)^{\otimes 2} \to \FHNA(V)
   \]
is defined as
   \[
   \mu_F\left((x_{1,n})_\tau, (x_{n+1, n+m})_\sigma\right) \,\buildrel \text{def} \over=\,
   (x_{1,n+m})_{\tau \vee \sigma}.
   \]
The linear map
   \[
   \alpha_F \colon \FHNA(V) \to \FHNA(V)
   \]
is defined by the rules:
   \begin{enumerate}
   \item $\alpha_F \vert V = \alpha_V$, and
   \item $\alpha_F\left((x_{1,n})_\tau\right) = (x_{1,n})_{\tau \lbrack 1 \rbrack}$ for $n \geq 2$.
   \end{enumerate}
Let $\iota \colon V \hookrightarrow \FHNA(V)$ denote the obvious inclusion map.

To show that $\FHNA$ is the left adjoint of $E$, let $(A, \mu_A, \alpha_A)$ be a Hom-nonassociative algebra, and let $f \colon V \to A$ be a morphism of Hom-modules.  We must show that there exists a unique morphism $g \colon \FHNA(V) \to A$ of Hom-nonassociative algebras such that $g \circ \iota = f$.  Define a map $g \colon \FHNA(V) \to A$ by
   \begin{equation}
   \label{eq:g}
   g\left((x_{1,n})_\tau\right) \,\buildrel \text{def} \over=\,
   \left(f(x_1) \cdots f(x_n)\right)_\tau,
   \end{equation}
where the right-hand side is defined as in \eqref{eq:product}.  It is clear that $g \circ \iota = f$.  Next we show that $g$ is a morphism of Hom-nonassociative algebras.

To show that $g$ commutes with $\alpha$, first note that $g \circ \alpha_F$ coincides with $\alpha_A \circ g$ when restricted to $V$, since $f$ commutes with $\alpha$.   For $n \geq 2$, we compute as follows:
   \[
   \begin{split}
   (g \circ \alpha_F)&((x_{1,n})_\tau)
   \,=\, g\left((x_{1,n})_{\tau \lbrack 1 \rbrack}\right)
   \,=\, \left(f(x_1) \cdots f(x_n)\right)_{\tau \lbrack 1 \rbrack} \\
   & \,=\, \alpha_A\left((f(x_1) \cdots f(x_n))_\tau\right)
   \,=\, (\alpha_A \circ g)((x_{1,n})_\tau).
   \end{split}
   \]
To show that $g$ is compatible with $\mu$, we compute as follows:
   \[
   \begin{split}
   (g \circ \mu_F)\left((x_{1,n})_\tau \otimes (y_{1,m})_\sigma\right)
   &= g\left((x_{1,n} \otimes y_{1,m})_{\tau \vee \sigma}\right) \\
   &= \left(f(x_1) \cdots f(x_n) f(y_1) \cdots f(y_m)\right)_{\tau \vee \sigma} \\
   &= \mu_A\left((f(x_1)\cdots f(x_n))_\tau \otimes (f(y_1) \cdots f(y_m))_\sigma\right) \\
   &= (\mu_A \circ g^{\otimes 2})\left((x_{1,n})_\tau \otimes (y_{1,m})_\sigma\right).
   \end{split}
   \]
This shows that $g$ is compatible with $\mu$ as well, so $g$ is a morphism of Hom-nonassociative algebras.

Note that, by a simple induction argument, every generator $(x_{1,n})_\tau \in \FHNA(V)$ can be obtained from $x_1, \ldots , x_n \in V$ by repeatedly taking products $\mu_F$ and applying $\alpha_F$.  Indeed, if $\tau = (\tau_1 \vee \tau_2) \lbrack r \rbrack \in \Twt_n$ as in \eqref{eq:tau}, then
   \[
   (x_{1,n})_{\tau} \,=\, \alpha_F^r\left((x_{1,p})_{\tau_1} (x_{p+1,p+q})_{\tau_2}\right).
   \]
The same argument then applies to $(x_{1,p})_{\tau_1}$ and $(x_{p+1,p+q})_{\tau_2}$, and so on.  This process has to stop after a finite number of steps, since in each step both $p$ and $q$ are strictly less than $n$.
Since $g$ is required to be a morphism of Hom-nonassociative algebras, this remark implies that $g$ is determined by its restriction to $V$, which must be equal to $f$.  This shows that $g$ is unique.
\end{proof}

We called $(\FHNA(V), \mu_F, \alpha_F) \in \HomNonAs$ the \emph{free Hom-nonassociative algebra} of $(V,\alpha_V) \in \HomMod$.  This object is the analogue in the Hom-nonassociative setting of the non-unital tensor algebra $T(V) = \oplus_{n \geq 1} V^{\otimes n}$.  Other free Hom-algebras can be obtained from the free Hom-nonassociative algebra.  One such example is given in Section \ref{subsec:freeHomAs}.

%%%%%%%%%%%%%%%%%%%%%%%%%%%%%%%%%%%%%%%%%%%%%%%%%
%%%%%%%%%%%%%%%%%%%%%%%%%%%%%%%%%%%%%%%%%%%%%%%%%
%%
%%     Enveloping algebras
%%
%%%%%%%%%%%%%%%%%%%%%%%%%%%%%%%%%%%%%%%%%%%%%%%%%
%%%%%%%%%%%%%%%%%%%%%%%%%%%%%%%%%%%%%%%%%%%%%%%%%
\section{Enveloping algebras of Hom-Lie algebras}
\label{sec:env}

The purpose of this section is to construct the enveloping Hom-associative algebra functor that is left adjoint to the functor $HLie$, which we first recall.

%%%%%%%%%%%%%%%%%%%%%%%%%%%%%%%
\subsection{The functor $HLie$}
\label{subsec:HLie}

A \emph{Hom-associative algebra} \cite[Definition 1.1]{ms} is a Hom-nonassociative algebra $(A, \mu, \alpha)$ such that the following $\alpha$-twisted associativity holds for $x, y, z \in A$:
   \begin{equation}
   \label{eq:alpha ass}
   \alpha(x) \left(y z\right) \,=\,
   \left(x y\right) \alpha(z).
   \end{equation}
As before, we abbreviate $\mu(x,y)$ to $xy$.  The full subcategory of $\HomNonAs$ whose objects are the Hom-associative algebras is denoted by $\HomAs$.

A \emph{Hom-Lie algebra} \cite[Definition 1.4]{ms} (first introduced in \cite{hls}) is a Hom-nonassociative algebra $(L, \lbrack -,-\rbrack, \alpha)$, satisfying the following two conditions:
   \begin{enumerate}
   \item $\lbrack x, y \rbrack = - \lbrack y, x \rbrack$ (skew-symmetry),
   \item $0 = \lbrack \alpha(x), \lbrack y, z \rbrack\rbrack + \lbrack \alpha(z), \lbrack x, y \rbrack \rbrack + \lbrack \alpha(y), \lbrack z, x \rbrack\rbrack$ (Hom-Jacobi identity)
   \end{enumerate}
for $x, y, z \in L$.  The full subcategory of $\HomNonAs$ whose objects are the Hom-Lie algebras is denoted by $\HomLie$.

Given a Hom-associative algebra $(A, \mu, \alpha)$, one can associate to it a Hom-Lie algebra $(HLie(A), \lbrack -,-\rbrack, \alpha)$ in which $HLie(A) = A$ as a $\mathbb{K}$-module and
   \[
   \lbrack x, y \rbrack \,\buildrel \text{def} \over=\, xy - yx
   \]
for $x, y \in A$ \cite[Proposition 1.7]{ms}.  The bracket defined is clearly skew-symmetric.  The Hom-Jacobi identity can be verified by writing out all $12$ terms and observing that their sum is $0$.

This construction gives a functor $HLie \colon \HomAs \to \HomLie$ that is the Hom-algebra analogue of the functor $Lie$ that associates a Lie algebra to an associative algebra via the commutator bracket.  The functor $Lie$ has as its left adjoint the enveloping algebra functor $U$.  We now construct the Hom-algebra analogue of the functor $U$, which is denoted by $\UHLie$.

%%%%%%%%%%%%%%%%%%%%%%%%%%%%%%%%
\subsection{Enveloping algebras}
\label{subsec:env}

Let $(L, \lbrack -,- \rbrack, \alpha)$ be a Hom-Lie algebra.  Consider the free Hom-nonassociative algebra $(\FHNA(L), \mu_F, \alpha_F)$ and the sequence of two-sided ideals, $I^1 \subset I^2 \subset \cdots \subset I^{\infty} \subset \FHNA(L)$, defined as follows.  Let $I^1$ be the two-sided ideal
   \[
   I^1 \,=\,
   \langle\, im(\mu_F \circ (\mu_F \otimes \alpha_F - \alpha_F \otimes \mu_F));~
   \lbrack x,y \rbrack - (xy - yx) \text{ for $x, y \in L$} \,\rangle.
   \]
Here the linear space $L$ is identified with its image under the inclusion $\iota \colon L \to \FHNA(L)$, and $xy$ denotes $\mu_F(\iota(x),\iota(y))$.  Inductively, we set
   \[
   I^{n+1} \,=\, \langle I^n \cup \alpha(I^n) \rangle, \quad
   I^\infty \,=\, \bigcup_{n \geq 1} I^n.
   \]

\begin{lemma}
\label{lem:env}
The submodule $I^\infty \subset \FHNA(L)$ is a two-sided ideal and is closed under $\alpha_F$.  The quotient $\FHNA(L)/I^\infty$, together with the induced maps of $\mu_F$ and $\alpha_F$, is a Hom-associative algebra.
\end{lemma}

\begin{proof}
Given elements $x \in I^\infty$ and $y \in \FHNA(L)$, we know that $x \in I^n$ for some $n < \infty$.  Therefore, both $xy$ and $yx$ lie in $I^n \subset I^\infty$.  The two-sided ideal $I^\infty$ is closed under $\alpha_F$ because, again, every element in $I^\infty$ must lie in some $I^n$, and $\alpha(I^n) \subset I^{n+1} \subset I^\infty$.

To show that $\FHNA(L)/I^\infty$, equipped with the induced maps of $\mu_F$ and $\alpha_F$, is a Hom-associative algebra, we must show that $\alpha$-associativity \eqref{eq:alpha ass} holds.  Let $x, y,$ and $z$ be elements in $L$.  Consider the diagram of $\mathbb{K}$-modules,
   \[
   \FHNA(L) \,\twoheadrightarrow\, \FHNA(L)/I^1 \,\twoheadrightarrow\, \FHNA(L)/I^\infty.
   \]
The first projection map sends the element $(\alpha(x)(yz) - (xy)\alpha(z)) \in \FHNA(L)$ to $0$, since it is in the image of the map $\mu_F \circ (\mu_F \otimes \alpha_F - \alpha_F \otimes \mu_F)$ and, therefore, in $I^1$.  It follows that the image of the element $(\alpha(x)(yz) - (xy)\alpha(z))$ in the quotient $\FHNA(L)/I^\infty$ is $0$ as well.  This shows that $\FHNA(L)/I^\infty$ is a Hom-associative algebra.
\end{proof}

From now on, we will denote the Hom-associative algebra $(\FHNA(L)/I^\infty, \mu_F, \alpha_F)$ by $(\UHLie(L), \mu, \alpha)$.  This defines a functor
   \[
   \UHLie \colon \HomLie \to \HomAs.
   \]

\begin{thm}
\label{thm:env}
The functor $\UHLie \colon \HomLie \to \HomAs$ is left adjoint to the functor $HLie$.
\end{thm}

\begin{proof}
Let $(L, \lbrack -,-\rbrack, \alpha_L)$ be a Hom-Lie algebra, and let $j \colon L \to \UHLie(L)$ be the composition of the maps
   \[
   L \,\buildrel \iota \over \hookrightarrow\, \FHNA(L) \,\buildrel pr \over \twoheadrightarrow\, \UHLie(L).
   \]
Let $(A, \mu_A, \alpha_A)$ be a Hom-associative algebra, and let $f \colon L \to HLie(A)$ be  a morphism of Hom-Lie algebras.  In other words, $f \colon L \to A$ is a linear map such that $f \circ \alpha_L = \alpha_A \circ f$ and
   \[
   f(\lbrack x, y \rbrack) \,=\, f(x)f(y) - f(y)f(x)
   \]
for $x, y \in L$.  We must show that there exists a unique morphism $h \colon \UHLie(L) \to A$ of Hom-associative algebras such that $f = h \circ j$ (as morphisms of $\mathbb{K}$-modules).

By Theorem \ref{thm:FreeHomNonAs}, there exists a unique morphism $g \colon \FHNA(L) \to A$ of Hom-nonassociative algebras such that $f = g \circ \iota$.  The map $g$ is defined in \eqref{eq:g}.  We claim that $g(I^\infty) = 0$.  It suffices to show that $g(I^n) = 0$ for all $n \geq 1$.  To see this, first note that $g(z) = 0$ for any element $z$ in the image of the map $\mu_F \circ (\mu_F \otimes \alpha_F - \alpha_F \otimes \mu_F)$, since $g$ commutes with both $\mu$ and $\alpha$ and $A$ satisfies $\alpha$-twisted associativity \eqref{eq:alpha ass}.  Moreover, for elements $x, y \in L$, we have that
   \[
   g(\lbrack x, y\rbrack - (xy - yx))
   \,=\, f(\lbrack x,y\rbrack) - (f(x)f(y) - f(y)f(x)) \,=\, 0.
   \]
It follows that $g(I^1) = 0$, again because $g$ commutes with $\mu$.  By induction, if $g(I^n) = 0$, then $g(\alpha(I^n)) = \alpha(g(I^n)) = 0$ as well.  Therefore, $g(I^{n+1}) = 0$, which finishes the induction step.  Since $g(I^n) = 0$ for all $n \geq 1$, it follows that $g(I^\infty) = 0$, as claimed.

The previous paragraph shows that $g$ factors through $\FHNA(L)/I^\infty = \UHLie(L)$.  In other words, there exists a linear map $h \colon \UHLie(L) \to A$ such that $g = h \circ pr$.  Since the operations $\mu$ and $\alpha$ on $\UHLie(L)$ are induced by the ones on $\FHNA(L)$, it follows that $h$ is also compatible with $\mu$ and $\alpha$.  In other words, $h$ is a morphism of Hom-associative algebras such that
   \[
   f \,=\, g \circ \iota \,=\, h \circ pr \circ \iota \,=\, h \circ j.
   \]
The uniqueness of $h$ follows exactly as in the last paragraph of the proof of Theorem \ref{thm:FreeHomNonAs}.  This finishes the proof of the Theorem.
\end{proof}

%%%%%%%%%%%%%%%%%%%%%%%%%%%%%%%%%%%%%%%%%%
\subsection{Free Hom-associative algebras}
\label{subsec:freeHomAs}

The construction of the functor $\UHLie$ can be slightly modified to obtain the free Hom-associative algebra functor.   Indeed, all we need to do is to redefine the ideals $I^n$ as follows.  Let $(V, \alpha)$ be a Hom-module.  Define $J^1 = \langle im(\mu_F \circ (\mu_F \otimes \alpha_F - \alpha_F \otimes \mu_F)) \rangle$, $J^{n+1} = \langle J^n \cup \alpha(J^n) \rangle$, and $J^\infty = \cup_{n\geq 1} J^n$.  Essentially the same argument as above shows that $J^\infty$ is a two-sided ideal that is closed under $\alpha$.  Moreover, the quotient module
   \begin{equation}
   \label{eq:FHAs def}
   \FHAs(V) \,\buildrel \text{def} \over=\, \FHNA(V)/J^\infty,
   \end{equation}
equipped with the induced maps of $\mu_F$ and $\alpha_F$, is the free Hom-associative algebra of $(V, \alpha)$.  In other words,
   \[
   \FHAs \colon \HomMod \to \HomAs
   \]
is the left adjoint of the forgetful functor $\HomAs \to \HomMod$.  The functor $\FHAs$ gives us a way to construct a Hom-associative algebra starting with just a Hom-module.

Conversely, if $(A, \mu, \alpha)$ is a Hom-associative algebra, then the adjoint of the identity map on $A$ is a surjective morphism $g \colon \FHAs(A) \twoheadrightarrow A$ of Hom-associative algebras.  The kernel of $g$ is a two-sided ideal in $\FHAs(A)$ that is closed under $\alpha$.  This allows us to write any given Hom-associative algebra $A$ as a quotient of a free Hom-associative algebra,
   \begin{equation}
   \label{eq:quotient}
   A \,\cong\, \FHAs(A)/\ker(g),
   \end{equation}
where the isomorphism is induced by $g$.

Other free Hom-algebras, including free Hom-dialgebras, free Hom-Lie algebras, and free Hom-Leibniz algebras, can be constructed similarly from the free Hom-nonassociative algebra.

%%%%%%%%%%%%%%%%%%%%%%%%%%%%%%%%%%%%%%%%%%%%%%%%%
%%%%%%%%%%%%%%%%%%%%%%%%%%%%%%%%%%%%%%%%%%%%%%%%%
%%
%%   Hom-dialgebras
%%
%%%%%%%%%%%%%%%%%%%%%%%%%%%%%%%%%%%%%%%%%%%%%%%%%
%%%%%%%%%%%%%%%%%%%%%%%%%%%%%%%%%%%%%%%%%%%%%%%%%
\section{Hom-dialgebras and Hom-Leibniz algebras}
\label{sec:HomDi}

The purposes of this section are (i) to introduce Hom-dialgebras and give some examples and (ii) to show how Hom-dialgebras give rise to Hom-Leibniz algebras.

%%%%%%%%%%%%%%%%%%%%%%%
\subsection{Dialgebras}
\label{subsec:dialg}

First we recall the definition of a dialgebra from \cite{loday}.  A \emph{dialgebra} $D$ is a $\mathbb{K}$-module equipped with two bilinear maps $\lprod$, $\rprod \colon D^{\otimes 2} \to D$, satisfying the following five axioms:
   \begin{equation}
   \label{eq:axioms}
   \begin{split}
   x \lprod (y \lprod z) &\,\buildrel (1) \over =\, (x \lprod y) \lprod z
   \,\buildrel (2) \over =\, x \lprod (y \rprod z), \quad 
   (x \rprod y) \lprod z \,\buildrel (3) \over =\, x \rprod (y \lprod z), \\
   (x \lprod y) \rprod z &\,\buildrel (4) \over =\, x \rprod (y \rprod z)
   \,\buildrel (5) \over =\, (x \rprod y) \rprod z
   \end{split}
   \end{equation}
for $x, y, z \in D$.  Many examples of dialgebras can be found in \cite[pp.\ 16-18]{loday}.

%%%%%%%%%%%%%%%%%%%%%%%%%%%
\subsection{Hom-dialgebras}
\label{susbec:HomDi}

We extend this notion to the Hom-algebra setting.  A \emph{Hom-dialgebra} is a tuple $(D, \lprod, \rprod, \alpha)$, where $D$ is a $\mathbb{K}$-module, $\lprod, \rprod \colon D^{\otimes 2} \to D$ are bilinear maps, and $\alpha \in \End(D)$, such that the following five axioms are satisfied for $x, y, z \in D$:
   \begin{equation}
   \label{eq:HomDiaxioms}
   \begin{split}
   \alpha(x) \lprod (y \lprod z) &\,\buildrel (1) \over =\, (x \lprod y) \lprod \alpha(z)
   \,\buildrel (2) \over =\, \alpha(x) \lprod (y \rprod z), \, 
   (x \rprod y) \lprod \alpha(z) \,\buildrel (3) \over =\, \alpha(x) \rprod (y \lprod z), \\
   (x \lprod y) \rprod \alpha(z) &\,\buildrel (4) \over =\, \alpha(x) \rprod (y \rprod z)
   \,\buildrel (5) \over =\, (x \rprod y) \rprod \alpha(z).
   \end{split}
   \end{equation}
We will often denote such a Hom-dialgebra by $D$.  A \emph{morphism} $f \colon D \to D^\prime$ of Hom-dialgebras is a linear map that is compatible with $\alpha$ and the products $\lprod$ and $\rprod$.  The category of Hom-dialgebras is denoted by $\HomDi$.

Note that if $D$ is a Hom-dialgebra, then, by axioms (1) and (5), respectively, both $(D, \lprod, \alpha)$ and $(D, \rprod, \alpha)$ are Hom-associative algebras.

%%%%%%%%%%%%%%%%%%%%%%%%%%%%%%%%%%%%%%%
\subsection{Examples of Hom-dialgebras}
\label{subsec:examples}

\begin{enumerate}
\item If $(A, \mu, \alpha)$ is a Hom-associative algebra, then $(A, \lprod, \rprod, \alpha)$ is a Hom-dialgebra in which $\lprod ~= \mu =~ \rprod$.
\item If $(D, \lprod, \rprod)$ is a dialgebra, then $(D, \lprod, \rprod, \alpha = \Id_D)$ is a Hom-dialgebra.
\item This example is an extension of \cite[Example 2.2(d)]{loday}.  First we need some definitions.  Let $(A, \mu_A, \alpha_A)$ be a Hom-associative algebra, and let $(M, \alpha_M)$ be a Hom-module.  A \emph{Hom-$A$-bimodule} structure on $(M, \alpha_M)$ consists of:
    \begin{enumerate}
    \item a left $A$-action, $A \otimes M \to M$ ($a \otimes m \mapsto am$), and
    \item a right $A$-action, $M \otimes A \to M$ ($m \otimes a \mapsto ma$)
    \end{enumerate}
such that the following three conditions hold for $x, y \in A$ and $m \in M$: $\alpha_A(x)(ym) = (xy)(\alpha_M(m))$, $(mx)(\alpha_A(y)) = (\alpha_M(m))(xy)$, and $\alpha_A(x)(my) = (xm)(\alpha_A(y))$.  A \emph{morphism} $f \colon M \to N$ of Hom-$A$-bimodules is a morphism $f \colon (M, \alpha_M) \to (N, \alpha_N)$ of Hom-modules such that $f(am) = af(m)$ and $f(ma) = f(m)a$ for $a \in A$ and $m \in M$.

For example, if $g \colon A \to B$ is a morphism of Hom-associative algebras, then $B$ becomes a Hom-$A$-bimodule via the actions, $ab = g(a)b$ and $ba = bg(a)$, for $a \in A$ and $b \in B$.  In particular, the identity map $\Id_A$ makes $A$ into a Hom-$A$-bimodule, and $g \colon A \to B$ becomes a morphism of Hom-$A$-bimodules.

Now let $(M, \alpha_M)$ be a Hom-$A$-bimodule, and let $f \colon M \to A$ be a morphism of Hom-$A$-bimodules.  Then the tuple $(M, \lprod, \rprod, \alpha_M)$ is a Hom-dialgebra in which 
$m_1 \lprod m_2 = m_1f(m_2)$ and $m_1 \rprod m_2 = f(m_1)m_2$ for $m_1, m_2 \in M$.  The five Hom-dialgebra axioms \eqref{eq:HomDiaxioms} are easy to check.  For example, given elements $m_1, m_2, m_3 \in M$, we have that
   \begin{eqnarray*}
   (m_1 \lprod m_2) \lprod \alpha_M(m_3)
   &=& (m_1f(m_2))(f(\alpha_M(m_3))) 
   ~=~ (m_1f(m_2))(\alpha_A(f(m_3))) \\
   &=& \alpha_M(m_1)(f(m_2)f(m_3))
   ~=~ \alpha_M(m_1)f\left(m_2f(m_3)\right) \\
   &=& \alpha_M(m_1) \lprod (m_2 \lprod m_3).
   \end{eqnarray*}
This shows (1) in \eqref{eq:HomDiaxioms}.  The other four axioms are checked similarly.
\end{enumerate}

%%%%%%%%%%%%%%%%%%%%%%%%%%%%%%%%%%%%%%%%%%%%%%%%%%%%%%%%
\subsection{From Hom-dialgebras to Hom-Leibniz algebras}
\label{subsec:from}

Recall from \cite[Definition 1.2]{ms} that a \emph{Hom-Leibniz algebra} is a triple $(L, \lbrack -,- \rbrack, \alpha)$, in which $L$ is a $\mathbb{K}$-module, $\alpha \in \End(L)$, and $\lbrack -,- \rbrack \colon L^{\otimes 2} \to L$ is a bilinear map, that satisfies the Hom-Leibniz identity,
   \begin{equation}
   \label{eq:HomLeibnizId}
   \lbrack \lbrack x, y \rbrack, \alpha(z) \rbrack  \,=\,
   \lbrack \lbrack x, z \rbrack, \alpha(y) \rbrack +
   \lbrack \alpha(x), \lbrack y, z \rbrack \rbrack,
   \end{equation}
for $x, y, z \in L$.  The full subcategory of $\HomNonAs$ whose objects are Hom-Leibniz algebras is denoted by $\HomLeib$.

Note that Hom-Lie algebras are examples of Hom-Leibniz algebras.  Also, if $\alpha = \Id_L$ in a Hom-Leibniz algebra $(L, \lbrack -,- \rbrack, \alpha)$, then $(L, \lbrack -,-\rbrack)$ is called a Leibniz algebra \cite{loday0,loday1,loday,lp,lp2}, which is a non-skew-symmetric version of a Lie algebra.  In \cite[Proposition 4.2]{loday}, Loday showed that a dialgebra gives rise to a Leibniz algebra via a version of the commutator bracket (see \eqref{eq:commutatorbracket} below).  The result below is the Hom-algebra analogue of that result.

\begin{prop}
\label{prop:HomDi2HomLeib}
Let $(D, \lprod, \rprod, \alpha)$ be a Hom-dialgebra.  Define a bilinear map $\lbrack -,-\rbrack \colon D^{\otimes 2} \to D$ by setting
   \begin{equation}
   \label{eq:commutatorbracket}
   \lbrack x, y \rbrack \,\buildrel \text{def} \over=\, x \lprod y - y \rprod x.
   \end{equation}
Then $(D, \lbrack -,-\rbrack, \alpha)$ is a Hom-Leibniz algebra.
\end{prop}

\begin{proof}
We write down all twelve terms involved in the Hom-Leibniz identity \eqref{eq:HomLeibnizId}:
   \[
   \begin{split}
   \lbrack \lbrack x,y \rbrack, \alpha(z)\rbrack
   &\,=\, (x \lprod y) \lprod \alpha(z) - (y \rprod x) \lprod \alpha(z) - \alpha(z) \rprod (x \lprod y) + \alpha(z) \rprod (y \rprod x), \\
   \lbrack \lbrack x,z \rbrack, \alpha(y) \rbrack
   &\,=\, (x \lprod z) \lprod \alpha(y) - (z \rprod x) \lprod \alpha(y) - \alpha(y) \rprod (x \lprod z) + \alpha(y) \rprod (z \rprod x), \\
   \lbrack \alpha(x), \lbrack y, z \rbrack \rbrack
   &\,=\, \alpha(x) \lprod (y \lprod z) - \alpha(x) \lprod (z \rprod y) - (y \lprod z) \rprod \alpha(x) + (z \rprod y) \rprod \alpha(x).
   \end{split}
   \]
Using the five Hom-dialgebra axioms \eqref{eq:HomDiaxioms}, it is immediate to see that \eqref{eq:HomLeibnizId} holds.
\end{proof}

We write $(HLeib(D), \lbrack -,-\rbrack, \alpha)$ for the Hom-Leibniz algebra $(D, \lbrack -,-\rbrack, \alpha)$ in Proposition \ref{prop:HomDi2HomLeib}.  This gives a functor
   \begin{equation}
   \label{eq:HLeib}
   HLeib \colon \HomDi \to \HomLeib,
   \end{equation}
which is the Hom-Leibniz analogue of the functor $HLie$ \cite[Proposition 1.7]{ms}.

%%%%%%%%%%%%%%%%%%%%%%%%%%%%%%%%%%%%%%%%%%%%%%%%%%%%%
%%%%%%%%%%%%%%%%%%%%%%%%%%%%%%%%%%%%%%%%%%%%%%%%%%%%%
%%
%%     Enveloping algebras
%%
%%%%%%%%%%%%%%%%%%%%%%%%%%%%%%%%%%%%%%%%%%%%%%%%%%%%%
%%%%%%%%%%%%%%%%%%%%%%%%%%%%%%%%%%%%%%%%%%%%%%%%%%%%%
\section{Enveloping algebras of Hom-Leibniz algebras}
\label{sec:env Leibniz}

The purpose of this section is to construct the left adjoint $\UHLeib$ of the functor $HLeib$. On the one hand, this is the Leibniz analogue of the functor $\UHLie$ (Theorem \ref{thm:env}).  On the other hand, this is the Hom-algebra analogue of the enveloping algebra functor of Leibniz algebras \cite{lp}.

As in the case of $\UHLie$, the construction of $\UHLeib$ depends on a suitable notion of trees, which we discuss next.

%%%%%%%%%%%%%%%%%%%%%%%%%%%%%
\subsection{Diweighted trees}
\label{subsec:diweighted}

By a \emph{diweighted $n$-tree}, we mean a pair $\tau = (\psi, w)$ in which:
   \begin{enumerate}
   \item $\psi \in T_n$ is an $n$-tree, called the \emph{underlying $n$-tree of $\tau$}, and
   \item $w$ is a function from the set of internal vertices of $\psi$ to the set $\mathbb{Z}_{\geq 0} \times \lbrace \lprod, \rprod \rbrace$.  We call $w$ the \emph{weight function of $\tau$}.
   \end{enumerate}
The set of diweighted $n$-trees is denoted by $\Tdi_n$.  As in the case of weighted trees, we have $\Tdi_1 = T_1$.  Every diweighted $n$-tree $\tau = (\psi, w)$ can be pictured by drawing the underlying $n$-tree $\psi$ and putting the weight $w(v)$ next to each internal vertex $v$ of $\psi$.

Let $\tau = (\psi, w) \in \Tdi_n$ and $\tau^\prime = (\psi^\prime, w^\prime) \in \Tdi_m$ be two diweighted trees.  Define the \emph{left grafting} to be the diweighted $(n+m)$-tree, $\tau \vee_l \tau^\prime \buildrel \text{def} \over= (\psi \vee \psi^\prime, \omega)$, where the weight function is given by
   \[
   \omega(v) \,=\,
   \begin{cases}
   w(v) & \text{if $v$ is an internal vertex of $\psi$}, \\
   w^\prime(v) & \text{if $v$ is an internal vertex of $\psi^\prime$}, \\
   (0, \lprod) & \text{if $v$ is the lowest internal vertex of $\psi \vee \psi^\prime$}.
   \end{cases}
   \]
The \emph{right grafting} $\tau \vee_r \tau^\prime$ is defined in exactly the same way, except that $\omega(v) = (0, \rprod)$ if $v$ is the lowest internal vertex of $\psi \vee \psi^\prime$.

Let $m$ be a non-negative integer.  Suppose that $\tau = (\psi, w)$ is a diweighted $n$-tree in which $w(v) = (s, \ast) \in \mathbb{Z}_{\geq 0} \times \lbrace \lprod, \rprod \rbrace$, where $v$ is the lowest internal vertex of $\psi$.  Define a new diweighted $n$-tree, $\tau\lbrack m \rbrack \buildrel\text{def}\over = (\psi, w\lbrack m \rbrack)$, in which the weight function is given by
   \[
   w\lbrack m \rbrack(u) \,=\,
   \begin{cases}
   w(u) & \text{if $u \not= v$}, \\
   (s + m, \ast) & \text{if $u = v$}.
   \end{cases}
   \]
In other words, $\tau\lbrack m \rbrack$ adds $m$ to the integer component of the weight of the lowest internal vertex of $\tau$.

Every diweighted $n$-tree $\tau = (\psi, w)$ can be written uniquely in the form
   \begin{equation}
   \label{eq:tau di}
   \tau \,=\, (\tau_1 \vee_\ast \tau_2) \lbrack m \rbrack,
   \end{equation}
where $\ast \in \lbrace l, r\rbrace$, $\tau_1 \in \Tdi_p$, $\tau_2 \in \Tdi_q$ with $p + q = n$, and $m$ is the integer component of the weight of the lowest internal vertex of $\tau$.  Every diweighted $n$-tree for $n \geq 2$ can be obtained from $n$ copies of the $1$-tree by iterating the operations $\vee_l$, $\vee_r$, and $\lbrack m\rbrack$ $(m \geq 0)$.

%%%%%%%%%%%%%%%%%%%%%%%%%%%%%%%%
\subsection{Enveloping algebras}
\label{subsec:env Leibniz}

Let $(V, \alpha_V)$ be a Hom-module.  Consider the module
   \begin{equation}
   \label{eq:Fbar}
   \Fbar(V) \,\buildrel \text{def} \over=\,
   \bigoplus_{n \,\geq\, 1} \bigoplus_{\tau \,\in\, \Tdi_n} V^{\otimes n}_\tau,
   \end{equation}
where $V^{\otimes n}_\tau$ is a copy of $V^{\otimes n}$.  A generator $(x_1 \otimes \cdots \otimes x_n)_\tau \in V^{\otimes n}_\tau$ will be abbreviated to $(x_{1,n})_\tau$.  Define two bilinear operations $\lprod, \rprod \colon \Fbar(V)^{\otimes 2} \to \Fbar(V)$ by setting
   \[
   (x_{1,n})_\tau \ast (x_{n+1,n+m})_\sigma
   \,=\,
   \begin{cases}
   (x_{1,n+m})_{\tau \vee_l \sigma} & \text{if $\ast =~ \lprod$}, \\
   (x_{1,n+m})_{\tau \vee_r \sigma} & \text{if $\ast =~ \rprod$}.
   \end{cases}
   \]
Define a linear map $\alpha_F \colon \Fbar(V) \to \Fbar(V)$ by the rules: (1) $\alpha_F \vert V = \alpha_V$, and (2) $\alpha_F\left((x_{1,n})_\tau\right) = (x_{1,n})_{\tau \lbrack 1 \rbrack}$ for $n \geq 2$.

Note that $\Fbar(V)$ is the analogue of $\FHNA(V)$ (Theorem \ref{thm:FreeHomNonAs}) with two bilinear operations.  In $\Fbar(V)$, the \emph{two-sided ideal} $\langle S \rangle$ generated by a non-empty subset $S$ is the smallest sub-$\mathbb{K}$-module of $\Fbar(V)$ containing $S$ that is closed under both $\lprod$ and $\rprod$ (but not necessarily $\alpha$). It can be constructed as the sub-$\mathbb{K}$-module of $\Fbar(V)$ spanned by all the parenthesized monomials $x_1 \ast \cdots \ast x_n$ with $n \geq 1$ and $\ast \in \lbrace \lprod, \rprod \rbrace$ such that at least one $x_j$ lies in $S$.

Let $(L, \lbrack -,-\rbrack, \alpha)$ be a Hom-Leibniz algebra.  Define an increasing sequence of two-sided ideals, $I^1 \subset I^2 \subset \cdots \subset I^\infty \subset \Fbar(L)$,  
as follows.  Set $I^1$ to be the two-sided ideal in $\Fbar(V)$ generated by the subset consisting of:
\[
\begin{split} 
&(1)\, im(\lprod \circ~ (\lprod \otimes \alpha) - \lprod \circ~ (\alpha \otimes \lprod)), \,
(2)\, im(\lprod \circ~ (\lprod \otimes \alpha) - \lprod \circ~ (\alpha \otimes \rprod)), \\
&(3)\, im(\lprod \circ~ (\rprod \otimes \alpha) - \rprod \circ~ (\alpha \otimes \lprod)), \,
(4)\, im(\rprod \circ~ (\alpha \otimes \rprod) - \rprod \circ~ (\rprod \otimes \alpha)), \\
&(5)\, im(\rprod \circ~ (\alpha \otimes \rprod) - \rprod \circ~ (\lprod \otimes \alpha)), \,
(6)\, \lbrack x, y \rbrack - (x \lprod y - y \rprod x)
\end{split}
\]
for $x, y \in L$.
In (6), $L$ is regarded as a submodule of $\Fbar(V)$ via the inclusion map $\iota \colon L \hookrightarrow \Fbar(V)$, and $x \ast y = \iota(x) \ast \iota(y)$ for $x, y \in L$ and $\ast \in \lbrace \lprod, \rprod \rbrace$.  The first five types of generators in $I^1$ correspond to the five Hom-dialgebra axioms \eqref{eq:HomDiaxioms}.  Inductively, set
   \begin{equation}
   \label{eq:Iinfty}
   I^{n+1} \,=\, \langle I^n \cup \alpha(I^n)\rangle, \quad
   I^\infty \,=\, \bigcup_{n \geq 1} I^n.
   \end{equation}

We are now ready for the Leibniz analogue of the enveloping Hom-associative algebra functor $\UHLie$.

\begin{thm}
\label{thm:envLeib}
Let $(L, \lbrack -,-\rbrack, \alpha)$ be a Hom-Leibniz algebra.  Then:
   \begin{enumerate}
   \item $I^\infty$ \eqref{eq:Iinfty} is a two-sided ideal in $\Fbar(L)$ that is closed under $\alpha$.
   \item The quotient module $\UHLeib(L) \,\buildrel \text{def} \over=\, \Fbar(L)/I^\infty$, 
   equipped with the induced maps of $\lprod$, $\rprod$, and $\alpha$, is a Hom-dialgebra.
   \item The functor $\UHLeib \colon \HomLeib \to \HomDi$ 
   is left adjoint to the functor $HLeib$ \eqref{eq:HLeib}.
   \end{enumerate}
\end{thm}

Since this Theorem can be proved by arguments that are essentially identical to those in Section \ref{sec:env}, we will omit the proof.

%%==============%%
%%==============%%
%%              %%
%%  References  %%
%%              %%
%%==============%%
%%==============%%
\sgsp


\begin{thebibliography}{99}

\bibitem{hls}J.T. Hartwig, D. Larsson, and S.D. Silvestrov, Deformations of Lie algebras using $\sigma$-derivations, J. Algebra \textbf{295} (2006), 314-361.

\bibitem{loday0}J.-L. Loday, Cyclic homology, Grundl. Math. Wiss. Bd. \textbf{301}, Springer, 1992.

\bibitem{loday1}J.-L. Loday, Une version non commutative des alg\`{e}bres de Lie: les alg\`{e}bres de Leibniz, Ens. Math. \textbf{39} (1993), 269-293.

\bibitem{loday}J.-L. Loday, Dialgebras, in:  Dialgebras and related operads, 7-66, Lecture Notes in Math. \textbf{1763}, Springer, 2001.

\bibitem{lp}J.-L. Loday and T. Pirashvili, Universal enveloping algebras of Leibniz algebras and (co)homology, Math. Ann. \textbf{296} (1993), 139-158.

\bibitem{lp2}J.-L. Loday and T. Pirashvili, The tensor category of linear maps and Leibniz algebras, Georgian Math. J. \textbf{5} (1998), 263-276.

\bibitem{lq}J.-L. Loday and D. Quillen, Cyclic homology and the Lie algebra homology of matrices, Comm. Math. Helv. \textbf{59} (1984), 565-591.

\bibitem{ms}A. Makhlouf and S. Silvestrov, Hom-algebra structures, J. Gen. Lie Theory Appl., to appear, \texttt{arXiv:math.RA/0609501}.

%\bibitem{myung}H.\ C.\ Myung, Lie-admissible algebras, Hadronic J.\ \textbf{1} (1978), 169-193.

%\bibitem{weibel}C.\ A.\ Weibel, An introduction to homological algebra, Cambridge studies in advanced mathematics \textbf{38}, Cambridge Univ.\ Press, Cambridge, UK, 1994.

\end{thebibliography}
\end{document}